\documentclass[a4paper,12pt]{article}
\usepackage[margin=1in]{geometry}  

\usepackage{graphicx}              
\usepackage{amsmath}
\usepackage{amscd}               
\usepackage{amsfonts} 
\usepackage{amssymb}             
\usepackage{amsthm}                

\newtheorem{theorem}{Theorem}[section]

\newtheorem{lemma}[theorem]{Lemma}
\newtheorem{proposition}[theorem]{Proposition}
\newtheorem{corollary}[theorem]{Corollary}

\newcommand{\cR}{\mathcal{R}} 

\newcommand{\R}{\mathbb{R}}      

\newcommand{\E}{\mathbb{E}}     

\newcommand{\cN}{\mathcal{N}}
\newcommand{\cH}{\mathcal{H}}
\newcommand{\cL}{\mathcal{L}}
\newcommand{\cI}{\mathcal{I}}
\newcommand{\cG}{\mathcal{G}}
\newcommand{\bs}[1]{{(#1)}}

\begin{document}
\title{Asymptotic independence of local eigenvalue distributions in macroscopically distinct regions}
\author{Yunjiang Jiang}
\maketitle
\section{Introduction}
	Consider the Wigner random matrix ensemble $X_N = (x_{ij})_{i,j=1}^N$ where
	\begin{enumerate}
	\item $x_{ij}$, $i < j$, are iid complex v alued with mean 0 and variance $1/N$.
	\item $x_{ii}$ are iid real-valued with mean 0 and variance $2/N$, and 
	\item $x_{ji} = \bar{x}_{ij}$, for $i < j$, i.e., the matrix is Hermitian.
	\item The distributions of $x_{ij}$ and $x_{ii}$ are supported on at least three points, and have finite fourth moments.
	\end{enumerate}
	
	The main result of the present paper is 
	\begin{theorem}\label{Main}
	
	Let $f_1, \ldots, f_p$ be compactly supported functions on $\R^{k_1}, \ldots, \R^{k_p}$ respectively, $-2 < a_1 < \ldots < a_p < 2$, and $I_i = [a_i -\epsilon_i, a_i + \epsilon_i]$ such that $I_i$'s are disjoint from each other. For a compactly function $g$ on $\R^k$, and $a \in (-2,2)$, denote
	\begin{align*}
	g[a](\lambda_1, \ldots, \lambda_N) = \sum_{1 \le i_1 < \ldots < i_k \le N} g(N \rho_{\rm{sc}}(a)(\lambda_{i_1}-a), \ldots, N \rho_{\rm{sc}}(a)(\lambda_{i_k}-a)).
	\end{align*}
	Then 
	\begin{align*}
	\lim_N \E  \prod_{i=1}^p \frac{1}{2\epsilon_i}\int_{u_i \in I_i}  f_i[u_i](\lambda_1, \ldots, \lambda_N)du_i - \prod_{i=1}^p \E  \frac{1}{2\epsilon_i}\int_{u_i \in I_i}  f_i[u_i](\lambda_1, \ldots, \lambda_N)du_i = 0.
	\end{align*}
	\end{theorem}
	Let $\Delta_i := \Delta_i^{(N)}= [E_i, E_i + \frac{c_i}{N}]$, $i=1,\ldots, p$ for any fixed $p$, where $E_i$ are distinct points in $(-2,2)$ and $c_i$ are fixed real numbers independent of $N$. Also let $\cN(\Delta)$ be the number of eigenvalues in $\Delta$. 
	 
	 As a corollary we have
	 \begin{corollary}
	 The following averaged asymptotic independence of eigenvalue counting functions hold (see the Theorem 3 of \cite{BiDeNa} for GUE)
	\begin{align*}
	 &\lim_{N \to \infty} \frac{1}{\prod_{i=1}^p (2 \epsilon_i)}\int_{E_1 \in I_1} \ldots \int_{E_p \in I_p} \mathbb{P}[\cN(\Delta_1) = l_1, \ldots, \cN(\Delta_p) = l_p]dE_1 \ldots dE_p \\
	 &- \prod_{i=1}^p \frac{1}{2\epsilon_i} \int_{E_i \in I_i}\mathbb{P}[\cN(\Delta_i) = l_i] dE_i = 0.
	 \end{align*}
	\end{corollary}
	
	Vast amount of efforts has been devoted to proving that the eigenvalue distributions of $X_N$ as $N \to \infty$ have the same limit as that of the Gaussian Unitary Ensemble (GUE), where $x_{ij}$ are iid of a particular distribution, namely Gaussians. Whenever a property P of the GUE limiting spectrum is also true for any Wigner ensemble defined above, we say that P is universal.
	\\
	For single eigenvalue marginals, one can prove that the empirical spectral measure $\frac{1}{N} \sum_{i=1}^N \delta_{\lambda_i}$ converges almost surely to the semi-circle distribution, whose density is given by $\rho_{sc}(x) = \frac{1}{\pi} \sqrt{4 -x^2} 1_{|x| \le 2}$. Thus the eigenvalues of $X_N$ are supported on the interval $[-2,2]$ in the limit.
	\\
	The next order of investigation is the joint distribution of more than one eigenvalues. In the case of Hermitian ensembles, one could label the eigenvalues from 1 to $N$ according to their relative orders on $\R$. This perspective is very useful in proving the universality of local eigenvalue distribution, to be defined more precisely below (see \cite{ERSY}).
	\\
	A more natural point of view is to count how many eigenvalues lie in a region $\Delta \subseteq [-2,2]$, the size of which could shrink as $N \to \infty$, in order to have a nondegenerate limiting distribution. More generally, we define the so-called m-point correlation function $\rho_m(x_1, \ldots, x_m) = \rho_m^{(N)}(x_1, \ldots, x_m)$ by 
	\begin{align*}
	\int_{\R^m} f(x_1, \ldots, x_m) \rho_m(x_1, \ldots, x_m) dx_1 \ldots dx_m = \E \sum_{1 \le i_1 < \ldots < i_m \le N} f(\lambda_{i_1}, \ldots, \lambda_{i_m}),
	\end{align*}
	where $f$ is any bounded continuous test function on $\R^m$. Notice that $\rho_m$ is not a probability density on $\R^m$ since it integrates to $\frac{N!}{(N-m)!}$ against the constant $1$ function.
	\\
	It turns out that for GUE as well as a large class of unitarily invariant Hermitian ensembles(but not Wigner matrices in general), these m-point correlation functions can be expressed as a determinant of an m-dimensional matrix, whose elements are given by a fixed kernel as follows:
	\begin{align*}
	\rho_m^{\rm{GUE}(N)} (x_1, \ldots, x_m) = \det(K_N(x_i,x_j))_{i,j=1}^m,
	\end{align*}
	where $K_N(x,y) = N \frac{\psi_N(x \sqrt{N}) \psi_{N-1}(Y \sqrt{N}) - \psi_{N-1}(x \sqrt{N}) \psi_N(y \sqrt{N})}{x -y}$, and $\psi_k = \frac{e^{-x^2/4}}{\sqrt{\sqrt{2\pi} k!}} h_k(x)$ is the Hermite polynomial $h_k$ normalized and multiplied by square root of the Gaussian weight. 
	\\
	Thus in order to study the limiting behavior of $\rho_m$ as $N \to \infty$< it suffices to understand $K_N(x,y)$. It turns out, using asymptotics from orthogonal polynomial theory, 
	\begin{align*}
	\lim_{N \to \infty} \frac{1}{\rho_{sc}(x) N} K_N(x,x+ \frac{c}{N}) = \frac{\sin c}{c}.
	\end{align*}
	In the past decade or so, through the various effort of Brezin and Hikami \cite{BrHi}, Johansson \cite{Jo01}, Erdos, Ramirez, Peche, Schlein, Yau, Yin \cite{ERSY}, \cite{EYY}, and Tao,Vu \cite{TaoVu}, the universality of the sine kernel limit in the bulk of the eigenvalue spectrum have been established for very general Wigner ensembles, with only a mild condition of subexponential decay for individual entry distribution. 
	\\
	In this article, we look at a different property of the eigenvalues, namely independence in the limit $N \to \infty$ of two or more eigenvalue clusters located at macroscopically separate positions on $(-2,2)$. According to our best knowledge, the only results in this drection so far are given by
	\begin{enumerate}
	\item Bianchi, Debbah, and Najim \cite{BiDeNa}, who showed that 
	\begin{align*}
	\lim_{N \to \infty} | \mathbb{P}_N[ \cN(\Delta_1) = l_1, \ldots, \cN(\Delta_p) = l_p
	- \prod_{i=1}^p \mathbb{P}_N[\cN(\Delta_i) = l_i]| = 0,
	\end{align*}
	where $\Delta_i = [ E_i, E_i + \frac{c_i}{N}]$ for $E_1 < E_2 \ldots < E_p \in (-2,2)$, i.e., distinct locations in the limit spectral support, and $c_i \in \R$. If $c_i$ are negative, then we interpret $\Delta_i = [E_i - \frac{c_i}{N}, E_i]$. Finally $\cN(\Delta)$ denotes the number of eigenvalues in $\Delta$, here under the GUE ensemble only. In fact they showed more: if $E_1 = -2$ and/or $E_p = 2$, then the same result above holds provided $\Delta_1 = [-2,-2 + \frac{c_1}{N^{1/3}}]$ and $\Delta_p = [2,2+ \frac{c_p}{N^{1/3}}]$, in accordance with the scaling at the spectral edge.
	\\
	In the appendix, we reprove their result in the bulk using the inclusion-exclusion principle, instead of their more sophisticated approach using Fredholm determinant and complex analytic tools.
	\item F. Boremann \cite{Bo10} reproved the asymptotic independence of the extreme eigenvalues $\lambda_1$ and $\lambda_N$ for GUE using abstract operator theoretic argument. He also obtained the next lowest order term in their joint distribution, and used it as a distinguishing statistics for GUE against other Wigner ensembles.
	\item Soshnikov's earlier result \cite{So99} on universality of the edge distribution of Wigner matrices whose entries have symmetric distribution and satisfy a sub-Gaussian moment assumption actually proves the asymptotic independence of $\lambda_1$ and $\lambda_N$ as well, using his high power moment method.
	
	\item Borodin, Okounkov, and Olshanski \cite{BoOkOl} established the asymptotic connection between the bulk of GUE spectrum and the Frobenius coordinates of a random partition distributed according to the Plancherel measure. They mentioned that the distribution of Frobenius coordinates at far away regions are asymptotically independent, which presumably is proved in their paper also.
	\end{enumerate}
	
	In the present article we show that the result from \cite{BiDeNa} is universal.
	
	\section{Outline of the proof}
	We borrow heavily ideas from \cite{ERSY}. We illustrate the proof with $p=2$ only. The more general case follows the same pattern, without additional assumptions, unlike the main results in \cite{EPRSY}.
	
	The first step is to show that the eigenvalues confined in two clusters $\Lambda_1 \cup \Lambda_2 := \{\lambda_{L_1+1}, \ldots, \lambda_{L_1 + n}\} \cup \{ \lambda_{L_2+1}, \ldots, \lambda_{L_2 + n}\}$ converges in total variation distance to their joint stationary measure conditioned on the positions of the remaining $\lambda_j$'s $\Lambda \setminus (\Lambda_1 \cup \Lambda_2)$, under the mean-reverting Dyson's Brownian motion in time $t << 1$. Here $n = N^\delta$ for some small $\delta \in (0,1)$, and $0 < \lim_N \frac{L_1}{N} < \lim_N \frac{L_2}{N} < 1$, so the two clusters are anchored at distinct points on $(-2,2)$ in the limit. 
	\\
	As in \cite{ERSY}, we prove convergence in total variations using relative entropy bound by establishing a log-sobolev inequality for the invariant measure, which can be achieved very elegantly using Bakry-Emery theory. We will prove that the Hessian of the corresponding local Hamiltonian has spectrum bounded below by a constant in the next section.
	\\
	Next we show that for almost all $L_1 < L_2$, the conditional equilibrium measure of $\Lambda_1 \cup \Lambda_2$ is close in total variation to a product measure on $\Lambda_1 \times \Lambda_2$, where the components are the invariant measure on $\Lambda_1$ conditioned on $\Lambda \setminus (\Lambda_1 \cup \Lambda_2)$ and that on $\Lambda_2$ conditioned on the same event, for asymptotically almost all $\Lambda \setminus (\Lambda_1 \cup \Lambda_2)$ configurations.
	\\
	 In \cite{ERSY}, it has been proved that the correlation kernel of $\Lambda_1$ (or $\Lambda_2$) when rescaled so that $\rm{supp} \Lambda_1 = [-1,1]$ converges locally to the sine-kernel. This relied on results from \cite{LeLu} where it was shown that sine-kernel is universal in a broad class of unitarily invariant ensembles whose matrix element joint density is given by $\exp(- N \rm{Tr} H(X))$. 

Notice that $L_1$ and $L_2$ are the ordinal index of the eigenvalues, which are not directly related to their rescaled positions in $[-2,2]$. However, a concentration result gives that the eigenvalues do not deviate from their semi-circle law predicted positions by more than $\mathcal{O}(N)$. More recently, \cite{EYY} has shown that the deviation cannot exceed a polylogarithmic amount with high probability uniformly over all eigenvalues up to the edge. So one can essentially treat $L_1,L_2$ the same as $E_1,E_2$ the corresponding energy levels. 

Finally under the framework of \cite{ERSY} as well as their sequels, an averaging over a macroscopic range of energy levels $E_1$ and $E_2$ is needed for the above sine-kernel convergence, because the convergence to local equilibrium cannot be shown for all $L_1$. 

Thus their result was of the form
\begin{align} \label{ERSY main}
&\lim_{\delta\to 0} \lim_{N \to \infty} \frac{1}{2\delta} \int_{E_0 - \delta}^{E_0 + \delta} dE \int_\R \int_\R da db O(a,b) \frac{1}{N^2\rho_{sc}^2(E)} \rho_2(E + \frac{a}{\rho_{sc}(E) N}, E + \frac{b}{\rho_{sc} N}) \\
&= \int_\R g(u) [1-(\frac{\sin \pi u}{\pi u})^2] du,
\end{align}
	where $O(a,b):= g(a-b) h(\frac{a+b}{2})$ and $g,h$ are both compactly supported with $\int h(x) dx = 1$. Note that $h$ is just a local bump function that does not appear on the right hand side. If we do have 
	\begin{align*}
	\frac{1}{N^2\rho_{sc}(E)^2} \rho_2(E + \frac{a}{\rho_{sc}(E) N}, E + \frac{b}{\rho_{sc} N}) = 1-(\frac{\sin \pi u}{\pi u})^2,
	\end{align*}
	then \eqref{ERSY main} would follow by a change of variable $(a,b) \mapsto (a-b,a+b)$.
	\\
	Using the very recent rigidity result of Erdos, Yau and Yin \cite{EYY}, one gets that the support of $\Lambda_1$, $[\lambda_{L_1},  \lambda_{L_1 +n+1}]$ contains a subinterval of $[-2,2]$ of size $N^{-\gamma}$ for some $\gamma < 1-\delta$ (recall $n = N^\delta$) with high probability. Hence provided we know $L_1$ is a good index, sine-kernel convergence without averaging can be established in the interval $[\cR^{-1}(L_1) + \frac{1}{3} N^{\delta-1}, \cR^{-1}(L_1+n+1) - \frac{1}{3}N^{\delta-1}]$.  
	\\ 
	Moreover, if we have two eigenvalue clusters $\Lambda_1$ and $\Lambda_2$ based at two good eigenvalue indices $L_1,L_2$, they would contain intervals $I_1 := [E_1, E_1 + N^{-\gamma}]$ and $I_2:= [E_2, E_2 + N^{-\gamma}]$ respectively with high probability. So conditioned on $\Lambda \setminus(\Lambda_1 \cup \Lambda_2)$, which stays away from $I_1 \cup I_2$ with high probability for the same reason, the distribution of the counting functions $1_{\cN(\Delta_1) = l_1}$ and $1_{\cN(\Delta_2) = l_2}$ with $\Delta_1 \subset I_1$ and $\Delta_2 \subset I_2$ are almost independent for almost all $\Lambda \setminus (\Lambda_1 \cup \Lambda_2)$ configurations, by local convergence to equilibrium. 
	\\
	This is however not yet enough to show $1_{\cN(\Delta_1) = l_1}$ and $1_{\cN(\Delta_2) = l_2}$ are independent unconditionally. For that one needs 
	\begin{align*}
	\mathbb{P}[ \cN(\Delta_1) = l_1, \cN(\Delta_2) = l_2| \Lambda \setminus (\Lambda_1 \cup \Lambda_2)] = f(l_1, l_2) + o(1)
	\end{align*}
	asymptotically almost surely, for a fixed function $f$. The relation between universality of$\cN(\Delta)$ and that of the local correlations, as established in the appendix, reduces the problem to showing that $\rho_m^{(N)}(x_1, \ldots, x_m | \Lambda \setminus (\Lambda_1 \cup \Lambda_2))$ for $x_i \in \Delta_1$ and for $x_i \in \Delta_2$ both converge to a constant function almost surely for each $m$. But this is essentially proved in \cite{ERSY} and strengthened by results in \cite{EYY} mentioned above. Finally we integrate over all configurations of $\Lambda \setminus(\Lambda_1 \cup \Lambda_2)$, using the boundedness of the integrand, to obtain
	\begin{align*}
	\E [ \mathbb{P}[ \cN(\Delta_1) = l_1, \cN(\Delta_2) = l_2| \Lambda \setminus(\Lambda_1 \cup \Lambda_2)]] = f(l_1,l_2) + o(1)
	\end{align*}
	This establishes the universality of asymptotic independence for ERSY-ensemble (i.e., Johansson matrix with vanishing GUE component). The crucial assumption for the above argument to work is the 3 conditions (2.4)-(2.6) imposed on the entry distribution of the Wigner component of the ERSY-ensemble. Those ensure that the global entropy is bounded by $N^{1+\epsilon}$ for $\epsilon > 1/4$ at the end of a burn-in period of $t = O(1/N)$. 
	\\
	The next step is to remove the Gaussian divisibility condition by invoking the powerful four-moment theorem of Tao and Vu, which states that if two Wigner ensembles $X_N$ and $Y_N$ have entries with the same first four moments, then their joint eigenvalue distributions for finitely many eigenvalues are very close in the sense of weak convergence.
	\begin{theorem} (\cite{TaoVu}) 
	Given that $X_N$ and $Y_N$ match up to fourth moments, and let $\lambda_i$ be the eigenvalues of $X_N$ and $\lambda'_i$ those of $Y_N$, then for some universal constant $K, c > 0$,
	\begin{align*}
	|\E G(\lambda_{i_1},\ldots, \lambda_{i_k}) - \E G(\lambda'_{i_1}, \ldots, \lambda'_{i_k})| < K n^{-c}
	\end{align*}
	for all $G$ with $\nabla^j G(x)| \le n^{1+c}$, where $j = 1,\ldots, 5$.
	\end{theorem}
	Note that the version above differs from the original theorem of Tao and Vu by the power of $n$ in the derivative bound of $G$ because they scale Wigner matrices by a factor of $N$.
	\\
	If we take $f(x_1, \ldots, x_{l_1})$ supported on $\Delta_1$ and $g(x_1, \ldots, x_{l_2})$ supported on $\Delta_2 > \Delta_1$, satisfying the derivative condition above, and denote
	\begin{align*}
	\E (fg) :=\frac{(N-l_1-l_2)!}{N!} \E \sum_{i_1 < \ldots < i_{l_1} < j_1 < \ldots < j_{l_2}} f(\lambda_{i_1}, \ldots , \lambda_{i_{l_1}}) g(\lambda_{j_1}, \ldots, \lambda_{j_{l_2}}), 
	\end{align*}
		then 
		\begin{align*}
		|\E_{X_N}(fg) - \E_{Y_N}(fg) | < K N^{-c}.
		\end{align*}
		From here it is easy to derive that the correlation function of the form \begin{align*}
		\frac{(N- l_1 - l_2)!}{N!}\rho^{Y_N}_{l_1+l_2}(x_1, \ldots, x_{l_1}, y_1, \ldots, y_{l_2}) - \frac{(N-l_1)!}{N!} \rho^{Y_N}_{l_1}(x_1, \ldots, x_{l_1})\frac{(N-l_2)!}{N!} \rho^{Y_N}_{l_2}(y_1, \ldots, y_{l_2})  = o(1),
		\end{align*}
		 for $x_i \in \Delta_1$ and  $y_j \in \Delta_2$, provided same holds for $X_N$. Thus if a matrix ensemble has the same first four moments as some Johansson matrix ensemble, it also enjoy the property of asymptotic independence. As pointed out in \cite{TaoVu}, if the entry distribution is supported on at least 3 points, then one can always find a Johansson matrix that match the first four moments, which concludes the proof. 
		
	\section{Convergence to local equilibrium for more than one eigenvalue clusters}
In this section, we prove convergence in total variation distance of the joint distribution of two mesoscopic eigenvalue clusters to the conditional equilbrium measure given the other eigenvalues. For ease of notation, we will stick with two eigenvalue clusters in this section. Our notation will be largely based on \cite{ERSY}, except that we will denote by $X_N^t := (1-t)^{1/2} W + t^{1/2} U$ the Johansson matrix ensemble. Let $L_1^{(N)}, L_2^{(N)} \in [N]$ be two sequences of integers such that $\lim_N  \frac{L_1^{(N)}}{N} = a$ and $\lim_N \frac{L_2^{(N)}}{N} =b$ both exist and 
\begin{align*}
 0 < a < b < 1.
\end{align*}
Let $n = N^\delta$ where $\delta \in (0,1)$ is fixed, but can be chosen as small as one likes. For simplicity, we will suppress the dependence on $N$ in the above notation whenever possible. Consider the following two contiguous eigenvalue clusters
\begin{align*}
 \Lambda_1 &:= (\lambda_{L_1 + 1}, \ldots, \lambda_{L_1 + n}) \\
 \Lambda_2 &:= (\lambda_{L_2 + 1}, \ldots, \lambda_{L_2 + n}).
\end{align*}
 Clearly these two clusters are disjoint for all large $N$. Define the relative entropy between two probability measures on the same space by 
\begin{align*}
 S(\mu,\nu) := \int_\Omega \log d\nu \frac{d\nu}{d\mu} = \int_\Omega d\mu f log f,
\end{align*}
whenever $f$ exists, where $f = \frac{d\nu}{d\mu}$ is the Radon-Nikodym derivative of $\nu$ with respect to $\mu$, and $S(\mu,\nu) = \infty$ otherwise.
	
Also define total variation distance between $\mu$ and $\nu$ by
\begin{align*}
 \| \mu - \nu\|_{\rm{TV}} = \sup_{A \subset \Omega} |\mu(A) - \nu(A)|. 
\end{align*}
When $f$ exists, one can alternatively define
\begin{align*}
 \| \mu - \nu\|_{\rm{TV}} = \int_\Omega |f -1| d\mu.
\end{align*}
Thus one expects total variation distance to be bounded by relative entropy. Indeed we have
\begin{align*}
 \| \mu - \nu\|_{\rm{TV}}^2 \le 2 S(\mu,\nu).
\end{align*}
For the reader's convenience, we include a proof of this fact in the appendix.

It  turns out that entropy analysis is well-suited for Dyson Brownian motion, via the so-called Bakry-Emery theory, which leads to a logarithmic sobolev inequality that relates the relative entropy with Dirichlet form applied to the square root of the running density of the Markov process (see \cite{GuZe} for an excellent introduction to the subject; also \cite{ERSY} section 5 for its application to DBM). 

In \cite{ERSY} three conditions (2.4) - (2.6) are assumed about the initial Wigner ensemble, which basically stipulates that the distribution of the individual matrix entries satisfies a logarithmic sobolev inequality and has sub-exponential decay, and that its fourier transform has polynomial decay for a sufficiently large polynomial degree.

These distributions are flexible enough that their first four moments can be anything that comes from a distribution with at least three points in the support (see \cite{TaoVu} Corollary 30). We will assume the same set of conditions in the next three sections. Under these conditions, we have the following bound on the global entropy after a burn-in period of $\mathcal{O}(1/N)$ of the eigenvalues of the Wigner copy $W$ evolving under the DBM:
\begin{lemma} 
(Lemma 5.1 of \cite{ERSY})
 Let $f_t$ be the density of the eigenvalue distribution of $X_N^t: =(1-t)^{1/2} W + t^{1/2} U$ with respect to the GUE spectral measure $\mu(d\lambda_1, \ldots, d\lambda_N) = Z_N^{-1} e^{-N \sum_{i=1}^N \lambda_i^2 - \sum_{i<j} \log (\lambda_i - \lambda_j)^2} d\lambda_1 \ldots d\lambda_N$, and where $U$ is a standard $GUE(N)$ matrix. Then for every $\alpha > 1/4$,
\begin{align*}
 S(f_{N^{-1}}\mu,\mu) \le C N^{1 + \alpha}.
\end{align*}
\end{lemma}
 
Next using (5.4) of \cite{ERSY}, we can bound the global Dirichlet form of $\sqrt{f_t}$ at time $t > N^{-1}$ by
\begin{align*}
 D(\sqrt{f_t}) \le \frac{S(f_{N^{-1}})}{t - N^{-1}}.
\end{align*}
  This easily gives 
\begin{align*}
 D(\sqrt{f_t}) \le C N^{2 + \alpha} \tau^{-1},
\end{align*}
for $t = \tau N^{-1}$ and $\tau \ge 2$ (see (6.1) of \cite{ERSY}).

The reason to work with Dirichlet forms instead of entropy directly is that the former has convenient additive property, which allows us to spread the global Dirichlet form into local ones. Let $\Lambda$ denote the vector of all eigenvalues of a Hermitian matrix ensemble, arranged in increasing order. The notation $\Lambda \setminus(\Lambda_1 \cup \lambda_2)$ thus has obvious meaning, i.e., the vector of the external eigenvalues arranged in increasing order, and so does $\Lambda \setminus \Lambda_1$. We will always use the letter $y$ to stand for a realization of $\Lambda \setminus(\Lambda_1 \cup \Lambda_2)$ and $x^{(1)}$ and $x^{(2)}$ to denote that of $\Lambda_1$ and $\Lambda_2$ respectively. Also we denote collectively $x= (x^{(1)},x^{(2)})$. Unlike in \cite{ERSY}, we index the components of $y$ consecutively. In analogy with the definition of localized Dirichlet form in section 6.2 of \cite{ERSY}, we define the following 2-point conditional local Dirichlet form by
\begin{align*}
 D_{(L_1,L_2),y} (f) =D_{(L_1,L_2),y} (f,f) := \int_{\Omega^{L_1,L_2}_y} \frac{1}{2N} |\nabla_x f|^2 d \mu_y^{L_1,L_2} (x),
\end{align*}
where $\mu_y^{L_1,L_2}$ is the invariant measure on $\Lambda_1 \cup \Lambda_2$ given $\Lambda \setminus (\Lambda_1 \cup \Lambda_2) =y$ and
\begin{align*}
 \Omega_y^{L_1,L_2} := ([y_{L_1}, y_{L_1 + 1}]^n)^{\uparrow} \times ([y_{L_2-n}, y_{L_2-n + 1}]^n)^{\uparrow},
\end{align*}
where for a subset $I \subseteq \R$, we denote by $(I^n)^{\uparrow}$ the set of increasing $n$-tuples with components in $I$.
Note that $[y_{L_2-n}, y_{L_2 -n+1}]$ is the confining range of $\Lambda_2$. Denote by $f_{t,y}(x)$ the conditional density of the eigenvalues $\Lambda_1 \cup \Lambda_2$ at time $t$ given the external configuration $y$, with respect to the conditional invariant distribution on $\Lambda_1 \cup \Lambda_2$ given $y$. By exactly the same computation as in \cite{ERSY} between (6.17) and (6.18), we have
\begin{align} \label{average local Dirichlet}
\E_{f_t} D_{(L_1,L_2),y}(\sqrt{f_{t,y}(x)}) = \frac{1}{8N} \sum_{j=1}^n \int dx dy \frac{|\nabla_{x_j} f_t(y,x)|^2 + |\nabla_{x_{j+n}} f_t(y,x)|^2}{f_t(y,x)} u(y,x).
\end{align}
where $u(y,x)$ is the GUE eigenvalue density with respect to the Lebesgue measure on the positive orthant in $\R^N$. The expectation above is taken with respect to the external configuration $y$. 
\\
Now we sum the above equation \eqref{average local Dirichlet} over $L_1 \in [Na-n_2, Na + n_2]$ and $L_2 \in [Nb-n_2, Nb + n_2]$ where $n_2 = N^{\delta_2}$ and $\delta < \delta_2 < 1$. This amounts to summing the Dirichlet form for each individual particle in the range at most $n$ times. Thus we get
\begin{align*}
 \frac{1}{4n_2^2} \sum_{L_1 \in [Na-n_2, Na + n_2]} \sum_{L_2 \in \in [Nb-n_2, Nb + n_2]} \E_{f_t} D_{(L_1,L_2),y} (\sqrt{f_{t,y}(x)}) &\le \frac{n^2}{4N^{2 \delta_2}} D(\sqrt{f_t})\\
&\le C \frac{n^2}{n_2^2} N^{2+\alpha} \tau^{-1} 
\end{align*}
for some constant C. Thus by the discrete Markov inequality, if we let
\begin{align*}
 \mathcal{G}_N(a,b) &:= \{(L_1,L_2) \in [Na-n_2,Na+n_2] \times [Nb - n_2, Nb + n_2]: \E_{f_t} D_{(L_1,L_2),y} (\sqrt{f_{t,y}(x)}) \\
&\le C N^{2 + \alpha} \frac{n^3}{n_2^2} \tau^{-1}\}
\end{align*}
be the set of good pairs of eigenvalue indices, then for some constant $c$,
\begin{align*}
 \frac{|\mathcal{G}_N(a,b)|}{4n_2^2} \ge 1 - \frac{c}{n}.
\end{align*}
Furthermore, if $Y_{L_1,L_2} := \{y \in (\R^{N-2n})^{\uparrow}: D_{(L_1,L_2),y} (\sqrt{f_{t,y}}) \le C N^{2 + \alpha} \frac{n^4}{n_2^2} \tau^{-1} \}$, then by another application of Markov inequality, for $(L_1,L_2) \in \mathcal{G}_N(a,b)$, 
\begin{align*}
 \mathbb{P}_{f_t} [ Y_{L_1,L_2}] \ge 1 - c n^{-2}.
\end{align*}
Next we compute the $y$-conditional local Hamiltonian $\cH_y(x)$, which defines the $y$-conditional local equilibrium measure $\mu_y(dx) = \exp(-\cH_y(x))dx$:
\begin{align} \label{local H}
 \cH_y(x) &= N[ \frac{1}{2} \sum_{i=1}^n [(x_i^{(1)})^2 + (x_i^{(2)})^2] - \frac{1}{N} \sum_{1 \le i < j \le n} \log (x_i^{(1)} - x_i^{(1)})^2 - \frac{1}{N} \sum_{1 \le i < j \le n} \sum_{1 \le i < j \le n} \log (x_i^{(2)} \\
&- x_j^{(2)})^2 - \frac{1}{N} \sum_{i \in [n]}\sum_{j \in [n]} \log (x_i^{(1)} - x_j^{(j)})^2 - \frac{1}{N} \sum_{k=1}^{N-2n} \sum_{i=1}^n [\log (x_i^{(1)} - y_k)^2 + \log (x_i^{(2)} - y_k)^2 ]].
\end{align}
We can compute the Hessian of $\cH_y$ as follows:
\begin{align*}
 \partial_{x_i^{(1)}} \cH_y(x) = N x_i^{(1)} - \sum_{j\in [n] \setminus \{i\}} \frac{1}{x_i^{(1)} - x_j^{(1)}} - \sum_{j=1}^n \frac{2}{x_i^{(1)} - x_j^{(2)}} - \sum_{k=1}^{N-2n} \frac{2}{x_i^{(1)} - y_k}.
\end{align*}
Hence
\begin{align*}
 \partial_{x_j^{(1)}} \partial_{x_i^{(1)}} \cH_y(x) &= -\delta_{i \neq j} \frac{1}{(x_i^{(1)} - x_j^{(1)})^2} + \delta_{i = j} ( N + \sum_{k \neq i} \frac{1}{(x_i^{(1)} - x_k^{(1)})^2} \\
&+ \sum_{k=1}^n \frac{2}{(x_i^{(1)} - x_j^{(2)})^2} + \sum_{k=1}^{N-2n} \frac{2}{(x_i^{(1)} - y_k)^2}),\\
\partial_{x_j^{(2)}} \partial_{x_i^{(1)}} \cH_y(x) &= -\frac{2}{(x_i^{(1)} - x_j^{(2)})^2}.
\end{align*}
The other mixed partials $\partial_{x_i^{(2)}} \partial_{x_j^{(1)}}$ and $\partial_{x_i^{(2)}} \partial_{x_j^{(2)}}$ are obtained similarly. To obtain a spectral lower bound for $\rm{Hess} \cH_y$, we take a test vector $v \in \R^{2n}$, 
\begin{align*}
&\langle v, \text{Hess} \cH_y(x) v \rangle = \sum_{i=1}^n [N + \sum_{k \neq i} \frac{1}{(x_i^{(1)} - x_k^{(1)})^2} + \sum_{k=1}^n \frac{2}{(x_i^{(1)} - x_j^{(2)})^2} + \sum_{k=1}^{N-2n} \frac{2}{(x_i^{(1)} - y_k)^2}]v_i^2\\
 &+ \sum_{i=1}^n [N + \sum_{k \neq i} \frac{1}{(x_i^{(2)} - x_k^{(2)})^2} + \sum_{k=1}^n \frac{2}{(x_i^{(2)} - x_j^{(1)})^2} + \sum_{k=1}^{N-2n} \frac{2}{(x_i^{(2)} - y_k)^2}] v_{i+n}^2\\
 &- \sum_{i=1}^n \sum_{j \neq i} [ \frac{v_i v_j}{(x_i^{(1)} -x_j^{(1)})^2} + \frac{v_{i+n} v_{j+n}}{(x_i^{(2)} - x_j^{(2)})^2}] - 4 \sum_{i,j=1}^n \frac{v_i v_{j+n}}{(x_i^{(1)} -x_j^{(2)})^2}\\
&= N \| v\|^2 + \sum_{i=1}^n [ \sum_{k=1}^{N-2n} (\frac{2v_i^2}{(x_i^{(1)} -y_k)^2} + \frac{2v_{i+n}^2}{(x_i^{(2)} -y_k)^2}) + \sum_{i=1}^n \sum_{j\neq i} [ \frac{(v_i -v_j)^2}{  (x_i^{(1)} -x_j^{(1)})^2} + \frac{(v_{i+n} -v_{j+n})^2}{  (x_i^{(2)} -x_j^{(2)})^2} ] \\
&+ 2 \sum_{i,j=1}^n \frac{(v_i - v_{j+n})^2}{ (x_i^{(1)} -x_j^{(1)})^2}. 
\end{align*}
Thus for $x = (x^{(1)},x^{(2)})$ with $x^{(1)} \in ([y_{L_1},y_{L_1 +1}]^n)^\uparrow$ and $x^{(2)} \in ([y_{L_2-n},y_{L_2-n +1}]^n)^\uparrow$,
\begin{align*}
 \rm{Hess} \cH_y(x) \ge N + \inf_{\lambda \in [y_{L_1}, y_{L_1 +1}] \cup [y_{L_2-n}, y_{L_2 -n +1}]} C \sum_{k=1}^{N-2n} \frac{1}{(\lambda - y_k)^2}.
\end{align*}
 By rigidity of the eigenvalues of the Wigner ensemble $X_N^t$, proved in \cite{EYY}, the set $\Omega_2(L_1,L_2):=\{y \in  (\R^{N-2n})^\uparrow: \max\{y_{L_1+1} - y_{L_1}, y_{L_2-n+1} - y_{L_2 -n}\} > \frac{2n}{N}\} $ has vanishing probability. Indeed $\mathbb{P}_{X_N} (\Omega_2(L_1,L_2)) = O(N^{-\epsilon})$, for some $\epsilon > 0$. From this we easily conclude
\begin{align*}
 \rm{Hess} \cH_y(x) \ge \frac{CN^2}{n^2},
\end{align*}
 for almost all $y$. 

Next we invoke the Bakry-Emery theory as in \cite{ERSY} on the generator of the DBM $\cL := \frac{1}{2N} (\Delta - \langle \nabla \cH_y, \nabla \rangle)$, which gives 
\begin{align*}
 \partial_t D_{(L_1,L_2),y} (\sqrt{f_t}) \le -\frac{CN^2}{n^2 N} D(\sqrt{f_t}),
\end{align*}
 where again $f_t = e^{-t\cL} f_0$ is the the global density of the DBM at time $t$, and $f_0$ is the density of the joint eigenvalue distribution of starting ensemble $W_N$. This further implies the following logarithmic Sobolev inequality for any probability density $g$
\begin{align*}
 D(\sqrt{g}) \ge \frac{CN}{n^2} S(g),
\end{align*}
 where we abbreviated $S(g) := S(g,\mu_y^{L_1,L_2}, \mu_y^{L_1,L_2})$. See Appendix B for a proof of this. Note that the logarithmic Sobolev inequality bears no direct relation with any dynamical process that converges to the measure $\mu_y^{L_1,L_2}$. In fact we don't want to treat it from a dynamical point of view since the environment $y$ is changing during the time interval $[0,t]$ under DBM.

	Thus since $D_{(L_1,L_2),y}(\sqrt{f_t}) \le CN^{2+\alpha} \frac{n^4}{n_2^2} \tau^{-1}$, , we get for $(L_1,L_2) \in \mathcal{G}_N(a,b)$, 
\begin{align*}
 S(f_t) \le C N^{1+\alpha} \frac{n^6}{n_2^2} \tau^{-1}
\end{align*}
 So if $\tau >> N^{1+\alpha} n^6 / n_2^2$, which is still $o(1)$ if we choose $n_2$ sufficiently large and $n$ sufficiently small, we get $S(f_{t,y} \mu_y^{L_1,L_2}, \mu_y^{L_1,L_2}) = o(1)$, and hence
\begin{align*}
\| f_{t,y} \mu_y^{L_1,L_2}- \mu_y^{L_1,L_2} \|_{\rm{TV}} = o(1).
\end{align*}

\section{Asymptotic conditional independence of distant eigenvalue clusters}
In this section, we show that the invariant measure $\mu_y^{L_1,L_2}$ asymptotically converges to the product measure $\mu_y^{L_1} \otimes \mu_y^{L_2}$ with high probability. More precisely,
\begin{proposition}
 For each $N$ there is a set $\Omega_N^{L_1,L_2} \subseteq ([-2,2]^{N-2n})^\uparrow$ with $\mathbb{P}[ \Lambda \setminus (\Lambda_1^{L_1} \cup \Lambda_2^{L_2}) \in \Omega_N] \to 1$ such that for $(L_1,L_2) \in \mathcal{G}_N(a,b)$ the set of good index pairs,
\begin{align*}
 \| \mu_y^{L_1,L_2} - \mu_{y \cup z^{(2)}}^{L_1} \otimes \mu_{y \cup z^{(1)}}^{L_2} \|_{\rm{TV}} < C N^{-\epsilon},
\end{align*}
for some $\epsilon > 0$ and universal constant $C$. Here $z^{(1)} \in ([y_{L_1},y_{L_1 +1}])^\uparrow$ and $z^{(2)} \in ([y_{L_2-n},y_{L_2 - n+1}])^\uparrow$ are arbitrary. 
\end{proposition}
\begin{proof}
 Recall the 2-point local equilibrium measure anchored at $L_1, L_2$ has the following $y$-conditional density
\begin{align*}
 \mu_y^{L_1,L_2}(dx) = Z_{N,(L_1,L_2),y}^{-1} \exp( -\cH_y(x)) dx,
\end{align*}
where $\cH_y(x):=\cH_y^{L_1,L_2}(x)$ is given by \eqref{local H}. The key observation is that $\cH_y^{L_1,L_2}(x)$ splits approximately into a sum of two decoupled Hamiltonians:
\begin{align*}
 \cH^{L_1,L_2}_y(x) = \cH^{L_1}_{y \cup z^{(2)}}(x) + \cH^{L_2}_{y \cup z^{(1)}}(x) + R_z(x)
\end{align*}
where 
\begin{align*}
 \cH^{L_1}_{y \cup z^{(2)}}(x) &= N[ \sum_{i=1}^n \frac{1}{2} x_i^2 - \frac{1}{N} \sum_{1 \le i < j \le n} \log (x_j -x_i)^2 - \frac{1}{N} \sum_{k=1}^{N-2n} \sum_{i=1}^n (x_i - y_k)^2\\
& - \frac{1}{N} \sum_{k=1}^n \sum_{i=1}^n \log (x_i - z_k^{(2)})^2] \\
 \cH^{L_2}_{y \cup z^{(1)}}(x) &= N[ \sum_{i=1}^n \frac{1}{2} x_{i+n}^2 - \frac{1}{N} \sum_{1 \le i < j \le n} \log (x_{j+n} -x_{i+n})^2 \\
&- \frac{1}{N} \sum_{k=1}^{N-2n} \sum_{i=1}^n (x_{i+n} - y_k)^2 - \frac{1}{N} \sum_{k=1}^n \sum_{i=1}^n \log (x_{i+n} - z_k^{(2)})^2],
\end{align*}
and
\begin{align*}
 R_z(x) = \sum_{k=1}^n \sum_{i=1}^n \log (x_i - z_k^{(2)})^2 + \log (x_{i+n} - z_k^{(1)})^2 - \sum_{i=1}^n \sum_{j=1}^n \log(x_i - x_{j+n})^2.
\end{align*}
First observation is that $R_z(x)$ has small fluctuation for almost all $y$, i.e.,
\begin{align} \label{small fluctuation}
 \sup_{x,z} R_z(x) - \inf_{x,z} R_z(x) \le \sup_{x,z} C \sum_{k=1}^n \sum_{i=1}^n \frac{n/N}{x_i - z_k^{(2)}} +  \frac{n/N}{x_{i+n} - z_k^{(1)}} \le \frac{C n^3}{N},
\end{align}
where we used the first derivative estimate on the function $\log (x_i - z_k)^2$, and the rigidity result from \cite{EYY} on the external eigenvalue vector $y$ to lower bound the separation $x_i - z_k^{(2)}$ etc. 
\\ 
Next let $Z_1,Z_2$ be the partition function for $\cH^{L_1}_{y \cup z^{(2)}}$ and $\cH^{L_2}_{y \cup z^{(1)}}$ respectively, i.e., 
\begin{align*}
 Z_1 = \int_{([y_{L_1},y_{L_1 +1}]^n)^\uparrow} \exp(-\cH^{L_1}_{y \cup z^{(2)}}(x)) dx.
\end{align*}
Then comparing the following two quantities 
\begin{align*}
 Z_1 Z_2 &= \int_{([y_{L_1},y_{L_1 +1}]^n)^\uparrow} \int_{([y_{L_2-n},y_{L_2 -n +1}]^n)^\uparrow} \exp(-\cH^{L_1}_{y \cup z^{(2)}}(x) -\cH^{L_2}_{y \cup z^{(1)}}(x) ) dx^{(1)} dx^{(2)} \\
Z &:= Z_{N,(L_1,L_2),y}\\
&= \int_{([y_{L_1},y_{L_1 +1}]^n)^\uparrow} \int_{([y_{L_2-n},y_{L_2 -n +1}]^n)^\uparrow} \exp(-\cH^{L_1}_{y \cup z^{(2)}}(x) -\cH^{L_2}_{y \cup z^{(1)}}(x) + R_z(x)) dx^{(1)} dx^{(2)},
\end{align*}
we find for a fixed $z$ vector and any $x$ that
\begin{align*}
 Z_1 Z_2 \exp(R_z(x)) \exp(-C \frac{n^3}{N}) \le Z \le Z_1 Z_2 \exp(R_z(x)) \exp(C \frac{n^3}{N}). 
\end{align*}
This implies that the ratio of the density of $\mu_y^{L_1,L_2}$ with $\mu_{y \cup z^{(2)}}^{L_1} \otimes \mu_{y \cup z^{(1)}}^{L_2}$ satisfies the following bound for almost all $z$,

\begin{align*}
 \exp(-\frac{Cn^3}{N}) \le \frac{Z_1^{-1} \exp(-\cH_1(x^{(1)})) Z_2^{-1} \exp(-\cH_2(x^{(2)}))}{Z^{-1} \exp(-\cH(x))} \le \exp(\frac{Cn^3}{N}).
\end{align*}
Using now the $\cL^1$ interpretation of total variation distance, we easily arrive at
\begin{align*}
 \| \mu_y^{L_1,L_2} - \mu_{y \cup z^{(2)}}^{L_1} \otimes \mu_{y \cup z^{(1)}}^{L_2} \|_{\rm{TV}} &\le \int |\exp(\frac{Cn^3}{N}) -1| \mu_y^{L_1,L_2}(dx) \\
&\le \exp(\frac{Cn^3}{N}) -1.
\end{align*}
\end{proof}

\section{Unconditional independence}
In this section we will fix a sequence of good pairs of indices $(L_1,L_2) \in \mathcal{G}_N(a,b)$. By applying the main result from \cite{LeLu}, \cite{ERSY} showed that the process on $x^{(1)}= (x^{(1)}_1 < \ldots < x^{(1)}_n)$ under the measure $\mu_y^{L_1}(dx) = \exp(-\cH^{L_1}_y(x)) dx$ has the bulk correlation that converges to the sine-kernel when appropriately scaled so that $y_{L_1} = -2$ and $y_{L_1+1} = 2$, for all $y$ in some set $Y_N^{L_1}$ that has asymptotic probability 1 under the spectral distribution of the ERSY ensemble $J_N^t = (1-t)^{1/2} W_N + t^{1/2} H_N$ defined before. Here $\cH_y^{L_1}(x)$ is defined as before by
\begin{align*}
\cH_y(x) = N[ \frac{1}{2} \sum_{i=1}^n x_i^2 - \frac{1}{N} \sum_{i < j} \log (x_i-x_j)^2  - \frac{1}{N} \sum_{k=1}^{N-n} \sum_{i=1}^n \log (x_i - y_k)^2].
\end{align*}
Similarly one get sine kernel for $x^{(2)}=(x^{(2)}_1, \ldots, x^{(2)}_n)$, for $y \in Y_N^{L_2}$ with asymptotic probability 1. 
 
 The reason is that the Hamiltonian above is essentially the trace of some analytic function of the matrix ensemble, hence falls in the category of unitarily invariant ensembles, whose local correlation has long been conjectured to converge to the sine kernel, and proved under various conditions. To be explicit, we define the rescaling map
 \begin{align*}
 T^{(1)}_y(x) &= \frac{4(x-y_{L_1})}{y_{L_1+1}- y_{L_1}}\\
 T^{(2)}_y(x) &= \frac{4(x-y_{L_2-n})}{y_{L_2-n+1}-y_{L_2-n}}.
 \end{align*}
 Furthermore let $T^{(1)}_y(x^{(1)}) = (T^{(1)}_y(x^{(1)}_1), \ldots, T^{(1)}_y(x^{(1)}_n))$ and similarly define $T^{(2)}_y(x^{(2)})$. Then we have
 \begin{lemma}
 There exist sets $Y_N^{L_1,L_2} \subset \R^{N-2n}$ with asymptotic probability $1$ such that $y_N \in Y_N^{(L_1,L_2)}$ implies sine-kernel universality of local correlation for both $T^{(1)}_y(x^{(1)})$ and $T^{(2)}_y(x^{(2)})$. More precisely, for a compactly supported continuous function $f$, and $a \in (-2,2)$,
 \begin{align*}
 &\E_{y_N} \sum_{1 \le i_1 < \ldots < i_k \le n} f(n \rho_{\rm{sc}}(a)(T_y(x^{(1)}_{i_1}) -a),\ldots, n \rho_{\rm{sc}}(a)(T_y(x^{(1)}_{i_k}) -a))\\
  &= \int_{\R^k} f(z_1, \ldots, z_k) \rho_{\rm{sc}}^{(k)}( n \rho_{\rm{sc}}(a)(T_y(x^{(1)}_{i_1}) -a),\ldots, n \rho_{\rm{sc}}(a)(T_y(x^{(1)}_{i_k}) -a)) dz_1 \ldots dz_k.
  \end{align*}
  Obviously the same result holds for $T_y^{(2)}(x^{(2)})$. And in fact they can hold simultaneously by taking the intersection.
  \end{lemma}
  
  \begin{proof}
  We will adopt the notation 
  \begin{align*}
  f^y_a(T^\bs{1}(x^\bs{1})) = \sum_{1 \le i_1 < \ldots < i_k \le n} f(n \rho_{\rm{sc}}(a)(T_y(x^{(1)}_{i_1}) -a),\ldots, n \rho_{\rm{sc}}(a)(T_y(x^{(1)}_{i_k}) -a)).
  \end{align*}
  
   Let $K^{(1)}_y:= \{z^{(1)}:(y, z^{(1)}) \in Y_N^{L_1}\}$ and $K^{(2)}_y:= \{z^{(2)}:(y, z^{(2)}) \in Y_N^{L_2}\}$. For any $\epsilon \in (0,1)$, define the following projections of $Y_N^{L_1}$ and $Y_N^{L_2}$ onto their common domain:
  \begin{align*}
  Y_{N,\epsilon}^{(1)} &:= \{y \in \R^{N-2n}: \mathbb{P}_y[K^{(1)}_y] \ge \epsilon\} \\
  Y_{N,\epsilon}^{(2)} &:= \{y \in \R^{N-2n}: \mathbb{P}_y[K^{(2)}_y] \ge \epsilon\}, 
  \end{align*}
  Then for each fixed $\epsilon$, $\lim_N \mathbb{P}[Y_{N,\epsilon}^{(1)}] = 1$, which implies there is a sequence $\epsilon(N) \uparrow 1$ such that $Y_N^{L_1,L_2}:= Y_{N,\epsilon(N)}^{(1)}$ has asymptotic probability $1$ as well. Therefore, for $y \in Y_N^{L_1,L_2}$,
  \begin{align*}
 &\E_y f_a^y(T^\bs{1}(x^\bs{1}))\\
 &=  \int_{\R^k} f(z_1, \ldots, z_k) \rho_{\rm{sc}}^{(k)} (n \rho_{\rm{sc}}(a)(T_y(x^{(1)}_{i_1}) -a),\ldots, n \rho_{\rm{sc}}(a)(T_y(x^{(1)}_{i_k}) -a)) + o(1) + O(1-\epsilon(N)).
 \end{align*}
  \end{proof}
  
  Thus combined with the result on conditional asymptotic independence, we have for $\alpha, \beta \in(-2,2)$,
  \begin{align*}
  &\lim_N \E_{y_N} f^y_\alpha(T^\bs{1}(x^\bs{1})) g^y_\beta(T^\bs(2)(x^\bs{2}))\\
  &= \int_{([-2,2]^k)^\uparrow} \rho^\bs{k}_{\rm{sc}}(u_1, \ldots, u_k) f(\frac{n}{\rho_{\rm{sc}}(\alpha)}(u_1 - \alpha), \ldots, \frac{n}{\rho_{\rm{sc}}(\alpha)}(u_k - \alpha)) du_1 \ldots du_k\\
  & \int_{([-2,2]^l)^\uparrow} \rho^\bs{l}_{\rm{sc}}(v_1, \ldots, v_l) f(\frac{n}{\rho_{\rm{sc}}(\beta)}(v_1 - \beta), \ldots, \frac{n}{\rho_{\rm{sc}}(\beta)}(v_l - \beta)) dv_1 \ldots dv_l, 
  \end{align*}
  for $y_N \in Y_N^{L_1,L_2}$. This implies unconditional asymptotic independence of the $f^y_\alpha(T^\bs{1}(x^\bs{1}))$ and $g^y_\beta(T^\bs{2}(x^\bs{2}))$. 
  
  To obtain Theorem~\ref{Main}, we need to transfer the independence from indices to the energy levels. We will say a pair of energy levels $E_1 < E_2$ is covered by a pair of indices $L_1 < L_2$ if the classical locations of $[\lambda_{L_1}, \lambda_{L_1+n+1}]$ and $[\lambda_{L_2},\lambda_{L_2+n+1}]$ covers $E_1$ and $E_2$ respectively. More precisely, let $\cR^{-1}$ be the inverse cumulative function of the semicircle law on $[-2,2]$, then $L_1$ covers $E_1$ if
  \begin{align*}
  \cR^{-1}(L_1/N) < E_1 < \cR^{-1}((L_1 + n +1)/ N)
  \end{align*}
  and similarly for $L_2$ covering $E_2$. Furthermore we say $(E_1,E_2)$ is $\delta$-coverd by $(L_1,L_2)$ if 
  \begin{align*}
  \cR^{-1}(L_1/N) + \delta \frac{n}{N} < E_1 < \cR^{-1}((L_1 + n +1)/N) - \delta\frac{n}{N},
  \end{align*}
  i.e., the covering is padded.
  
  We will denote by $\cI_\delta(E)$ the set of all indices that $\delta$-cover $E$, and similarly denote $\cI_\delta(E_1,E_2)$ the set of all pairs of indices that $\delta$-cover $E_1,E_2$ respectively. Among $\cI_\delta(E_1,E_2)$ we denote by $\cG_\delta(E_1,E_2)$ the set of good pairs of indices. Now we have
  \begin{lemma}
  Let $I_i= [a_i - \epsilon_i,a_i + \epsilon_i]$, $i=1,2$ (for simplicity), be as in the statement of theorem~\ref{Main}. Denote by $K_{\epsilon,\delta} = \{(E_1,E_2) \in I_1 \times I_2: \frac{\cG_\delta(E_1,E_2)}{\cI_\delta(E_1,E_2)} > \epsilon\}$. Then 
  \begin{align*}
  \lim_N \frac{\lambda(K_{\epsilon,\delta})}{\lambda(I_1 \times I_2) }=1
  \end{align*}
  for every $\epsilon,\delta \in (0,1)$. Here $\lambda$ denotes Lebesgue measure on $\R$. 
  \end{lemma}
  This is proved essentially the same way as the previous lemma. So again we can let $\epsilon(N)$ be an increasing sequence going to $1$, and $\delta(N) \downarrow 0$, with 
  \begin{align*}
  \lim_N \frac{\lambda(K_{\epsilon(N),\delta(N)})}{\lambda(I_1 \times I_2)} = 1.
  \end{align*}
  Now to prove Theorem~\ref{Main}, it suffices to show that for any $\delta > 0$, and for asymptotically $\lambda$-almost all $(E_1,E_2) \in I_1 \times I_2$, the following is true:
  
  for almost all $(L_1,L_2) \in \cI_\delta(E_1,E_2)$, 
  \begin{align*}
  \E[ f_{\alpha_1}^y(T_y^\bs{1}(x^\bs{1})) g_{\alpha_2}^y(T_y^\bs{2}(x^\bs{2}))|y_{L_1},y_{L_1 +1}, y_{L_2 -n}, y_{L_2 -n+1}] = \E   f[a_1](\lambda_1, \ldots, \lambda_N) g[a_2](\lambda_1, \ldots, \lambda_N)
  \end{align*}
  where $\alpha_1 = \frac{4(a_1 - y_{L_1})}{y_{L_1 +1}- y_{L_1}}$, and $\alpha_2 = \frac{4(a_2 - y_{L_2 - n})}{y_{L_2 - n+1} - y_{L_2 -n}}$. We can then let $\delta$ go to $0$. But this follows directly from the previous lemma.


\begin{thebibliography}{99}
	
	\bibitem[BiDeNa]{BiDeNa} P. Bianchi, M. Debbah, J. Najim. Asymptotic independence in the spectrum of the gaussian unitary ensemble. ECP 2010. pp. 376-395.
	
 \bibitem[Bo10]{Bo10} F. Bornemann. Asymptotic independence of the extreme eigenvalues of Gaussian unitary ensemble. Journal of Mathematical Physics, 2010.	
	
	\bibitem[BoOkOl]{BoOkOl} A. Borodin, A. Okounkov, G. Olshanski. Asymptotics of Plancherel measures for symmetric groups. Journal of American Mathematical Society. Volume 13, Number 3, Pages 481-515.
	
	\bibitem[BrHi]{BrHi} E. Brezin, S. Hikami. Universal correlations for deterministic plus random Hamiltonians. Physical Review E, 1995. 
	
	\bibitem[EPRSY]{EPRSY} L. Erdos, S. Peche, J. A. Ramirez, B. Schlein, and H-T Yau. Bulk Universality for Wigner Matrices. Communications on Pure and Applied Mathematics, Vol. LXIII. 2010.
	
    \bibitem[ERSY]{ERSY} L. Erdos, J. A. Ramirez, B. Schlein, and H-T Yau. Universality of sine-kernel for Wigner matrices with a small Gaussian perturbation. EJP 2010. Vol. 15 no.18. 
    
    \bibitem[EYY]{EYY} L. Erdos, H-T Yau, J. Yin. Rigidity of Eigenvalues of Generalized Wigner Matrices. 2010. Arxiv 1007.4652v3.
     
\bibitem[GuZe]{GuZe} A. Guionnet, B. Zegarlinski. Lectures on logarithmic Sobolev inequalities. Seminaire de Probabilites: XXXVI, 2004.

\bibitem[Jo01]{Jo01} K. Johansson. Universality of the local spacing distribution in certain ensembles of Hermitian Wigner matrices.Communications in Mathematical Physics, 2001.

    \bibitem[LeLu]{LeLu} E. Levin, S.D.Lubinsky. Universality limits in the bulk for varying measures. Adv. Math. 219 (2008), 743-779.
    
    \bibitem[So99]{So99} A. Soshnikov. Universality at the Edge of the Spectrum of the Wigner random matrices.
    
    \bibitem[TaoVu]{TaoVu} T. Tao, V. Vu. Random matrices: universality of local eigenvalue statistics. To appear in Acta Math., Preprint. arXiv:0906.0510.
    
\end{thebibliography}
\end{document}